\documentclass[12pt]{article}

\usepackage{amsmath, epsfig, cite}
\usepackage{amssymb}
\usepackage{amsfonts}
\usepackage{latexsym}
\usepackage{float}

\newtheorem{thm}{Theorem}[section]

\newtheorem{conj}[thm]{Conjecture}
\newtheorem{lem}[thm]{Lemma}

\numberwithin{equation}{section}

\newcommand{\qed}{{\hfill$\square$}\medskip}

\begin{document}

%\linenumbers

\begin{center}
{\Large\bf On two conjectural supercongruences of\\[5pt]
Apagodu and Zeilberger}
\end{center}

\vskip 2mm \centerline{Ji-Cai Liu}
\begin{center}
{\footnotesize Department of Mathematics, Shanghai Key Laboratory of
PMMP, East China Normal University,\\ 500 Dongchuan Road, Shanghai
200241,
 People's Republic of China\\
{\tt jc2051@163.com} }
\end{center}

%%date: November 27, 2014
%%date : December 4, 2014

\vskip 0.7cm \noindent{\bf Abstract.}
Let the numbers $\alpha_n,\beta_n$ and $\gamma_n$ denote
\begin{align*}
\alpha_n=\sum_{k=0}^{n-1}{2k\choose k},\quad
\beta_n=\sum_{k=0}^{n-1}{2k\choose k}\frac{1}{k+1}\quad\text{and}\quad
\gamma_n=\sum_{k=0}^{n-1}{2k\choose k}\frac{3k+2}{k+1},
\end{align*}
respectively. We prove that for any prime $p\ge 5$ and positive integer $n$
\begin{align*}
\alpha_{np}&\equiv \left(\frac{p}{3}\right) \alpha_n \pmod{p^2},\\
\beta_{np}&\equiv \begin{cases}
\displaystyle \beta_n \pmod{p^2},\quad &\text{if $p\equiv 1\pmod{3}$},\\
-\gamma_n \pmod{p^2},\quad &\text{if $p\equiv 2\pmod{3}$},
\end{cases}
\end{align*}
where $\left(\frac{\cdot}{p}\right)$ denotes the Legendre symbol.
These two supercongruences were recently conjectured by Apagodu and Zeilberger.

\vskip 2mm
\noindent{\it MR Subject Classifications}: 11A07; 05A10

\section{Introduction}
In combinatorics, the $n$-th Catalan number and the $n$-th central binomial coefficient are given by ${2n\choose n}\frac{1}{n+1}$ and ${2n\choose n}$, respectively. Sun and Tauraso \cite[(1.7) and (1.9)]{ST2} have proved that for any prime $p\ge 5$
\begin{align}
&\sum_{k=0}^{p-1}{2k\choose k}\equiv \left(\frac{p}{3}\right)\pmod{p^2},\label{aa1}\\
&\sum_{k=0}^{p-1}{2k\choose k}\frac{1}{k+1}\equiv\frac{3}{2}\left(\frac{p}{3}\right)- \frac{1}{2}\pmod{p^2},\label{aa2}
\end{align}
where $\left(\frac{\cdot}{p}\right)$ denotes the Legendre symbol.

\noindent Recently, Apagodu and Zeilberger \cite{AZ} discussed some applications of the algorithm \cite{CHZ} for finding and proving some congruences (modulo primes) of indefinite sums of many combinatorial sequences. Let the numbers $\alpha_n,\beta_n$ and $\gamma_n$ denote
\begin{align*}
\alpha_n=\sum_{k=0}^{n-1}{2k\choose k},\quad
\beta_n=\sum_{k=0}^{n-1}{2k\choose k}\frac{1}{k+1}\quad\text{and}\quad
\gamma_n=\sum_{k=0}^{n-1}{2k\choose k}\frac{3k+2}{k+1},
\end{align*}
respectively. They have showed that for any prime $p\ge 5$ and positive integer $n$
\begin{align*}
\alpha_{np}&\equiv \left(\frac{p}{3}\right) \alpha_n \pmod{p},\\
\beta_{np}&\equiv \begin{cases}
\displaystyle \beta_n \pmod{p},\quad &\text{if $p\equiv 1\pmod{3}$},\\
-\gamma_n \pmod{p},\quad &\text{if $p\equiv 2\pmod{3}$}.
\end{cases}
\end{align*}

\noindent Apagodu and Zeilberger \cite[Super-Conjecture 1' and Super-Conjecture 2']{AZ} also conjectured supercongruence extensions of the above congruences,
for which the algorithm \cite{CHZ} was not applicable.

\begin{conj}\label{conj} (Apagodu-Zeilberger)
Suppose $p\ge 5$ is a prime and $n$ is a positive integer. Then
\begin{align}
\alpha_{np}&\equiv \left(\frac{p}{3}\right) \alpha_n \pmod{p^2},\label{aa3}\\
\beta_{np}&\equiv \begin{cases}
\displaystyle \beta_n \pmod{p^2},\quad &\text{if $p\equiv 1\pmod{3}$},\\
-\gamma_n \pmod{p^2},\quad &\text{if $p\equiv 2\pmod{3}$}.\label{aa4}
\end{cases}
\end{align}
\end{conj}

\noindent It is clear that \eqref{aa3} and \eqref{aa4} reduce to \eqref{aa1} and \eqref{aa2} when $n=1$.
The aim of this paper is to prove Conjecture \ref{conj}.
\begin{thm}\label{tt}
The supercongruences \eqref{aa3} and \eqref{aa4} are true.
\end{thm}

\noindent In the next section, we first prove some important lemmas. The proof of Theorem \ref{tt} will be given in Section 3. Throughout the paper, we use the convention that ${a\choose b}=0$ if $b>a$.
\section{Some lemmas}
\begin{lem}
Suppose $p\ge 5$ is a prime, $0\le k \le p-1$ and $r$ are non-negative integers. Then
\begin{align}
{2rp+2k\choose rp+k}\equiv {2r\choose r}{2k\choose k}(1+2rp(H_{2k}-H_k))\pmod{p^2},\label{bb1}
\end{align}
where the harmonic numbers are given by $H_n=\sum_{k=1}^n\frac{1}{k}$.
\end{lem}
{\noindent \it Proof.}
Note that
\begin{align}
{2rp+2k\choose rp+k}&={2rp\choose rp}\prod_{i=1}^{2k}(2rp+i)/\prod_{i=1}^{k}(rp+i)^2\notag\\
&\equiv {2r\choose r}\prod_{i=1}^{2k}(2rp+i)/\prod_{i=1}^{k}(rp+i)^2 \pmod{p^2},\label{bb2}
\end{align}
where we have utilized the Babbage's theorem \cite{Bab}
\begin{align}
{ap\choose bp}\equiv {a\choose b} \pmod{p^2}.\label{babb}
\end{align}

\noindent Now we consider the following rational function:
\begin{align}
f(x)=\prod_{i=1}^{2k}(2rx+i)/\prod_{i=1}^{k}(rx+i)^2.\label{bb3}
\end{align}
Taking the logarithmic derivative on both sides of \eqref{bb3} gives
\begin{align}
\frac{f'(x)}{f(x)}=\sum_{i=1}^{2k}\frac{2r}{2rx+i}-2\sum_{i=1}^{k}\frac{r}{rx+i}.
\label{bb4}
\end{align}
From \eqref{bb3} and \eqref{bb4}, we have
\begin{align*}
f(0)={2k\choose k}\quad \text{and}\quad f'(0)=2r{2k\choose k}(H_{2k}-H_{k}).
\end{align*}

\noindent Now we get the first two terms of the Taylor expansion for $f(x)$:
\begin{align}
f(x)={2k\choose k}+2rx{2k\choose k}(H_{2k}-H_{k})+\mathcal{O}\left(x^2\right).\label{bb5}
\end{align}
Combining \eqref{bb2} and \eqref{bb5}, we immediately get
\begin{align*}
{2rp+2k\choose rp+k}\equiv {2r\choose r}{2k\choose k}(1+2rp(H_{2k}-H_k))\pmod{p^2}.
\end{align*}
This completes the proof.
\qed

\begin{lem}
Let $n$ be an odd positive integer. Then
\begin{align}
&\sum_{k=0}^{n-1}{2k\choose k}\left({n-1\choose 2k}-(-1)^k{n-1\choose k}\right)=0,\label{bb6}\\
&\sum_{k=0}^{n-1}\frac{{2k\choose k}}{k+1}\left({n-1\choose 2k}-(-1)^k{n-1\choose k}\right)=\frac{1}{n+1}\sum_{k=0}^{n}(-1)^k{n+1\choose k}{2n-2k\choose n-k}.\label{bb7}
\end{align}
\end{lem}
{\noindent \it Proof.} These two combinatorial identities can be proved by Zeilberger algorithm \cite{PWZ}.
\qed

\begin{lem}
If $p\ge 5$ is a prime, then
\begin{align}
\sum_{k=0}^{p-1}{2k\choose k}\left(H_{2k}-H_{k}\right)\equiv 0 \pmod{p}.\label{bb8}
\end{align}
\end{lem}
{\noindent \it Proof.}
Note that for $0\le k\le p-1$
\begin{align}
(-1)^k{2k\choose k}{p-1\choose k}&=(-1)^k{2k\choose k}\frac{(p-1)(p-2)\cdots (p-k)}{k!}\notag\\
&\equiv {2k\choose k}(1-pH_k) \pmod{p^2}.\label{bb8-1}
\end{align}
Since ${2k\choose k}\equiv 0\pmod{p}$ for $(p+1)/2\le k \le p-1$, we have for $0\le k\le p-1$
\begin{align}
{2k\choose k}{p-1\choose 2k}&={2k\choose k}\frac{(p-1)(p-2)\cdots (p-2k)}{(2k)!}\notag\\
&\equiv {2k\choose k}(1-pH_{2k}) \pmod{p^2}.\label{bb8-2}
\end{align}
Letting $n=p$ in \eqref{bb6} and then using \eqref{bb8-1} and \eqref{bb8-2}, we get
\begin{align*}
p\sum_{k=0}^{p-1}{2k\choose k}\left(H_{2k}-H_k\right)\equiv 0 \pmod{p^2},
\end{align*}
which implies \eqref{bb8}.
\qed

\begin{lem}
If $p\ge 5$ is a prime, then
\begin{align}
\sum_{k=0}^{p-2}\frac{{2k\choose k}}{k+1}\left(H_{2k}-H_k\right)\equiv \frac{5}{2}+\frac{3}{2}\left(\frac{p}{3}\right) \pmod{p}.\label{bb9}
\end{align}
\end{lem}
{\noindent \it Proof.}
Similar to the proof of \eqref{bb8-1}and \eqref{bb8-2}, we have for $0\le k \le p-2$
\begin{align*}
\frac{(-1)^k}{k+1}{2k\choose k}{p-1\choose k}\equiv \frac{1}{k+1}{2k\choose k}(1-pH_k) \pmod{p^2},\\
\frac{1}{k+1}{2k\choose k}{p-1\choose 2k}\equiv \frac{1}{k+1}{2k\choose k}(1-pH_{2k}) \pmod{p^2}.
\end{align*}
Letting $n=p$ in \eqref{bb7}, and then applying the above two congruences and letting $k\to p-k$, we obtain
\begin{align}
p\sum_{k=0}^{p-2}\frac{{2k\choose k}}{k+1}\left(H_{k}-H_{2k}\right)&\equiv\frac{{2p-2\choose p-1}}{p}+\frac{1}{p+1}
\sum_{k=0}^{p}(-1)^k{p+1\choose k}{2p-2k\choose p-k}\notag\\
&=\frac{{2p-2\choose p-1}}{p}-\frac{1}{p+1}
\sum_{k=0}^{p}(-1)^k{p+1\choose k+1}{2k\choose k} \pmod{p^2}.\label{bb10}
\end{align}

\noindent Using \eqref{babb}, we obtain
\begin{align}
\frac{{2p-2\choose p-1}}{p}&=\frac{{2p\choose p}}{4p-2}
\equiv \frac{1}{2p-1}
\equiv -1-2p \pmod{p^2}.\label{bb10-1}
\end{align}

\noindent On the other hand,
\begin{align}
&\frac{1}{p+1}\sum_{k=0}^{p}(-1)^k{p+1\choose k+1}{2k\choose k}\notag\\
&=1-\frac{{2p\choose p}}{p+1}+\frac{1}{p+1}\sum_{k=1}^{p-1}(-1)^k{p+1\choose k+1}{2k\choose k}\notag\\
&\equiv -1+2p+p\sum_{k=1}^{p-1}\frac{(-1)^k}{k(k+1)}{p-1\choose k-1}{2k\choose k}\pmod{p^2}.\label{bb11}
\end{align}
Combining \eqref{bb10}, \eqref{bb10-1} and \eqref{bb11}, we obtain
\begin{align*}
\sum_{k=0}^{p-2}\frac{{2k\choose k}}{k+1}\left(H_{2k}-H_{k}\right)&\equiv 4+\sum_{k=1}^{p-1}\frac{(-1)^k}{k(k+1)}{p-1\choose k-1}{2k\choose k}\\
&\equiv 4-\sum_{k=1}^{p-1}\frac{1}{k(k+1)}{2k\choose k}\\
&=3-\sum_{k=1}^{p-1}\frac{1}{k}{2k\choose k}+\sum_{k=0}^{p-1}\frac{1}{k+1}{2k\choose k}\pmod{p},
\end{align*}
where we have utilized the fact that ${p-1\choose k-1}\equiv (-1)^{k-1}\pmod{p}$ in the second step. Then the proof of \eqref{bb9} follows from \eqref{aa2} and the following congruence of Pan and Sun \cite{PS} (see also \cite[Theorem 1.3]{ST1}):
\begin{align}
\sum_{k=1}^{p-1}\frac{1}{k}{2k\choose k}\equiv 0 \pmod{p}.\label{bb12}
\end{align}
\qed

\section{Proof of Theorem \ref{tt}}

{\noindent \it Proof of \eqref{aa3}.}
We first prove that for any non-negative integer $r$
\begin{align}
\sum_{k=rp}^{(r+1)p-1}{2k\choose k}\equiv \left(\frac{p}{3}\right){2r\choose r} \pmod{p^2}.\label{cc1}
\end{align}
Letting $k\to k+rp$ on the left-hand side of \eqref{cc1} and then using \eqref{bb1}, we have
\begin{align*}
\sum_{k=rp}^{(r+1)p-1}{2k\choose k}&=\sum_{k=0}^{p-1}{2rp+2k\choose rp+k}\\
&\equiv {2r\choose r}\sum_{k=0}^{p-1}{2k\choose k}(1+rp(2H_{2k}-2H_k))\pmod{p^2}.
\end{align*}
It follows from \eqref{aa1} and \eqref{bb8} that
\begin{align*}
\sum_{k=rp}^{(r+1)p-1}{2k\choose k}\equiv \left(\frac{p}{3}\right){2r\choose r} \pmod{p^2}.
\end{align*}
This concludes the proof of \eqref{cc1}. Taking the sum over $r$ from $0$ to $n-1$ on both sides of \eqref{cc1} gives
\begin{align*}
\sum_{k=0}^{np-1}{2k\choose k}\equiv \left(\frac{p}{3}\right)\sum_{k=0}^{n-1}{2k\choose k} \pmod{p^2},
\end{align*}
which is \eqref{aa3}.
\qed

{\noindent \it Proof of \eqref{aa4}.}
Letting $k\to k+rp$ and by \eqref{bb1}, we have
\begin{align}
\sum_{k=rp}^{(r+1)p-1}\frac{{2k\choose k}}{k+1}&=\sum_{k=0}^{p-1}\frac{{2rp+2k\choose rp+k}}{rp+k+1}\notag\\
&\equiv {2r\choose r}\sum_{k=0}^{p-2}
\frac{{2k\choose k}}{rp+k+1}\left(1+rp(2H_{2k}-2H_k)\right)
+\frac{{2(r+1)p-2\choose (r+1)p-1}}{(r+1)p}\pmod{p^2}.\label{cc1-1}
\end{align}
Note that
\begin{align}
\frac{1}{rp+k+1}\equiv \frac{1}{k+1}-\frac{rp}{(k+1)^2}\pmod{p^2}.\label{cc1-2}
\end{align}
By \eqref{babb}, we have
\begin{align}
\frac{{2(r+1)p-2\choose (r+1)p-1}}{(r+1)p}=\frac{{2(r+1)p\choose (r+1)p}}{4(r+1)p-2}\equiv-(2r+1){2r\choose r}\left(\frac{1}{r+1}+2p\right)\pmod{p^2}.\label{cc1-3}
\end{align}
Substituting \eqref{cc1-2} and \eqref{cc1-3} into \eqref{cc1-1} and then simplifying, we obtain
\begin{align}
\sum_{k=rp}^{(r+1)p-1}\frac{{2k\choose k}}{k+1}
&\equiv {2r\choose r}\sum_{k=0}^{p-2}
{2k\choose k}\left(\frac{1}{k+1}-\frac{rp}{(k+1)^2}+\frac{rp(2H_{2k}-2H_k)}{k+1}\right)\notag\\
&-(2r+1){2r\choose r}\left(\frac{1}{r+1}+2p\right)\pmod{p^2}.\label{cc2}
\end{align}

\noindent Using \eqref{aa2} and \eqref{bb10-1}, we have
\begin{align}
\sum_{k=0}^{p-2}\frac{{2k\choose k}}{k+1}&=\sum_{k=0}^{p-1}\frac{{2k\choose k}}{k+1}-\frac{{2p-2\choose p-1}}{p}\notag\\
&\equiv \frac{3}{2}\left(\frac{p}{3}\right)+\frac{1}{2}+2p \pmod{p^2}.\label{cc3}
\end{align}

\noindent By \eqref{bb12}, we obtain
\begin{align*}
\sum_{k=0}^{p-2}\frac{{2k\choose k}}{(k+1)^2}&=\sum_{k=1}^{p-1}{2k\choose k}\left(\frac{1}{2k-1}-\frac{1}{2k}\right)\\
&\equiv \sum_{k=1}^{p-1}{2k\choose k}\frac{1}{2k-1}\pmod{p}.
\end{align*}
Since ${2k\choose k}/(2k-1)=2{2k-2\choose k-1}/k$, by \eqref{aa2} and \eqref{bb10-1} we have

\begin{align}
\sum_{k=0}^{p-2}\frac{{2k\choose k}}{(k+1)^2}\equiv 2\sum_{k=0}^{p-2}\frac{{2k\choose k}}{k+1}
= 2\sum_{k=0}^{p-1}\frac{{2k\choose k}}{k+1}-\frac{2{2p-2\choose p-1}}{p}
\equiv 1+3\left(\frac{p}{3}\right) \pmod{p}.\label{cc4}
\end{align}

\noindent Combining \eqref{bb9} and \eqref{cc2}-\eqref{cc4}, we immediately get
\begin{align*}
\sum_{k=rp}^{(r+1)p-1}\frac{{2k\choose k}}{k+1}
&\equiv {2r\choose r}\left(\frac{3}{2}\left(\frac{p}{3}\right)-\frac{1}{2}-\frac{r}{r+1}\right)\pmod{p^2}.
\end{align*}

\noindent If $p\equiv 1\pmod{3}$, then $\left(\frac{p}{3}\right)=1$, and so
\begin{align}
\sum_{k=rp}^{(r+1)p-1}\frac{{2k\choose k}}{k+1}
&\equiv {2r\choose r}\frac{1}{r+1}\pmod{p^2}.\label{cc5}
\end{align}

\noindent If $p\equiv 2\pmod{3}$, then $\left(\frac{p}{3}\right)=-1$, and hence
\begin{align}
\sum_{k=rp}^{(r+1)p-1}\frac{{2k\choose k}}{k+1}
&\equiv -{2r\choose r}\frac{3r+2}{r+1}\pmod{p^2}.\label{cc6}
\end{align}
Taking the sum over $r$ from $0$ to $n-1$ on both sides of \eqref{cc5} and \eqref{cc6}, respectively, we obtain \eqref{aa4}. This completes the proof.
\qed

\end{document}